\documentclass{amsart}
\usepackage{amsmath,amssymb,epsfig,amscd}

\begin{document}

\def\R{{\mathbb R}}
\def\Z{{\mathbb Z}}
\def\Zt{{\mathbb Z}_2}
\def\Rn{{\mathbb R}^n}
\def\Rm{{\mathbb R}^m}
\def\U{{\mathcal U}}
\def\V{{\mathcal V}}
\def\tU{\widehat{\mathcal U}}
\def\tV{\widehat{\mathcal V}}
\def\fp{\text{\it fib-}\pi_0}
\def\acts{\hspace{.1cm}{
	\setlength{\unitlength}{.25mm}
	\linethickness{.08mm}
	\begin{picture}(8,8)(0,0)
	\qbezier(7,6)(4.5,8.3)(2,7)
	\qbezier(2,7)(-1.5,4)(2,1)
	\qbezier(2,1)(4.5,-.3)(7,2)
	\qbezier(7,6)(6.1,7.5)(6.8,9)
	\qbezier(7,6)(5,6.1)(4.2,4.4)
	\end{picture}\hspace{.1cm}
	}}
\def\acted{\hspace{.1cm}{
	\setlength{\unitlength}{.25mm}
	\linethickness{.08mm}
	\begin{picture}(8,8)(0,0)
	\qbezier(1,6)(3.5,8.3)(6,7)
	\qbezier(6,7)(9.5,4)(6,1)
	\qbezier(6,1)(3.5,-.3)(1,2)
	\qbezier(1,6)(1.9,7.5)(1.2,9)
	\qbezier(1,6)(3,6.1)(3.8,4.4)
	\end{picture}\hspace{.1cm}
	}}

\newtheorem{lemma}{Lemma}[section]
\newtheorem{xx}[lemma]{Example}
\newenvironment{example}{\begin{xx}\rm}{\end{xx}}
\newtheorem{prop}[lemma]{Proposition}
\newtheorem{df}[lemma]{Definition}
\newtheorem{idf}[lemma]{(Idea of) Definition}
\newtheorem{remark}[lemma]{Remark}

\title[Orbifolds from the Point of View of the Borel Construction]
{Orbispaces and Orbifolds from the Point of View of the Borel Construction,
a New Definition.}
\author{Andr\'e Henriques}

\maketitle
\vspace{-.8cm}

\section {Abstract}

``An orbifold is a space which is locally modeled on the quotient of a 
vector space by a finite group.'' This sentence is so easily said or
written that more than one person has missed some of the subtleties hidden by orbifolds.

Orbifolds were first introduced by Satake under the name ``V-manifold''
(cf \cite{Sat56} and \cite{Sat57}) and
rediscovered by Thurston who called them ``orbifolds''
(cf \cite{Thu}). Both of them used only faithful actions
to define their orbifolds
(these are the so called {\it reduced} orbifolds). This intuitive restriction
is very unnatural mathematically if we want for example
to study suborbifolds.

Orbispaces are to topological spaces
what orbifolds are to manifolds. They have been defined by
Haefliger (cf \cite{Hae84} and \cite {Hae90})
via topological groupoids (see also Kontsevich \cite{Kon95}).
In his papers, Haefliger defines the homotopy and (co)homology groups of
an orbispace to be those of the {\it classifying space} of the
topological groupoid. However, his definition of morphisms is 
complicated and technical. Ruan and Chen (cf \cite{Rua} and \cite{Che})
have tried to reformulate the definitions, but they don't
seem treat the non-reduced case in a satisfactory way.

Our idea is to define orbispaces and
orbifolds directly via their classifying spaces (called {\it Borel spaces}
in this paper, to stress the analogy with the Borel construction).
This approach allows us
to easily define morphisms of orbispaces and orbifolds with all the desired
good properties.

\section {Motivations}

Recall that, for a group $G$ acting on a space $X$,
the equivariant cohomology of $X$ is defined to be
$$
H^*_G(X)=H^*\big((X\times EG)/G\big),
$$
where $EG$ is a contractible space with free $G$-action
(see for example \cite{Gil}, chapt. 1).
The space $(X\times EG)/G$ is traditionally called the {\it Borel space}
of the $G$-space $X$. Now let $[X/G]$ be the orbispace defined by 
modding out $X$ by $G$; it is generally accepted that one should have
$H^*([X/G])=H^*_G(X)$. While the underlying space $X/G$ remembers
the geometry of $[X/G]$, the space $(X\times EG)/G$ is needed to
compute the algebraic topological invariants of $[X/G]$. We thus need
both of them for a good understanding of $[X/G]$, which motivates
the following definition.

\begin{idf}

An \underline{orbispace} $M$ is a surjective map
$p_M:PM\to QM$ such that for all $x\in QM$, $p_M^{-1}(x)$ is a $K(G,1)$
\footnote{
This definition is not complete (compare with definitions {\bf \ref{orbidef}}
and {\bf \ref{orbifolddef}})
because we don't have enough conditions on the topologies of $PM$ and $QM$.
For example, in the case where the homotopy type of $p_M^{-1}(x)$ does
not vary with $x$, we would also like $p_M$ to be a fibration.
}.
In other words, $\pi_k\big(p_M^{-1}(x)\big)=0$ for all $k\ge 2$.

An \underline{orbifold} $M$ is an orbispace such that $QM$ is covered
by open sets $\U$ with $p_M^{-1}(\U)\to\U$ isomorphic to
$(V\times EG)/G\to V/G$, where $G\acts V$ is a linear representation 
of a finite group.
\end{idf}

Note that the class of orbispaces naturally includes the class of spaces by
$X\mapsto\big(Id_X:X\to X\big)$.

\begin{example}[global quotient]\label{gl}
If $G\acts X$ is a quasi-free action (i.e. with discrete stabilizers)
of a compact Lie group on a manifold, then
$$
p:(X\times EG)/G\;\longrightarrow\;X/G
$$
is an orbifold, which will be denoted $[X/G]$.
\end{example}

Indeed, one can
check that for $x\in X/G$, $p^{-1}(x)=(Gx\times EG)/G=(\{x\}\times EG)/G_x
\simeq K(G_x,1)$, where $G_x$ is the stabilizer of $x$.

The way we defined orbispaces and orbifolds makes it obvious to what
orbispace maps $f:M\to N$
should be, namely commutative diagrams of the form
$$
\begin{CD}
PM	@>Pf>>	PN\\
@VVp_MV	@VVp_NV\\
QM	@>Qf>>	QN,
\end{CD}
$$
where $Qf$ is the underlying continuous map and $Pf$ is some lifting
of it to the respective Borel spaces. However, this provides {\it too many}
maps, because the Borel spaces are usually very big. It is thus
convenient to say that two maps $f$ and $g$ are {\it essentially the
same} if $Qf=Qg$ and if $Pf$ and $Pg$ are {\it homotopic over} $QN$.
More precisely, this means that one can deform
$Pf$ into $Pg$ by a homotopy that fixes the composition with $p_N$.
In other words, the homotopy changes the maps only in the
direction of the fibers of $p_N$.

Note that a map of orbifolds cannot be determined by
its underlying continuous map.
To illustrate why this is nessessary, consider the following paradox.

\begin{example}\label{paradoxe}
Let $M$ be the orbifold $S^1\times[\mathbb R/\mathbb Z_2]$ and let
$f$ and $g$ be the maps from $S^1$ to $M$ ``defined by'' $f(z)=(z,1)$ and
$g(z)=\big(z,Re(\sqrt{z})\big)$.
\vspace{.6cm}

\centerline{\epsfig{file=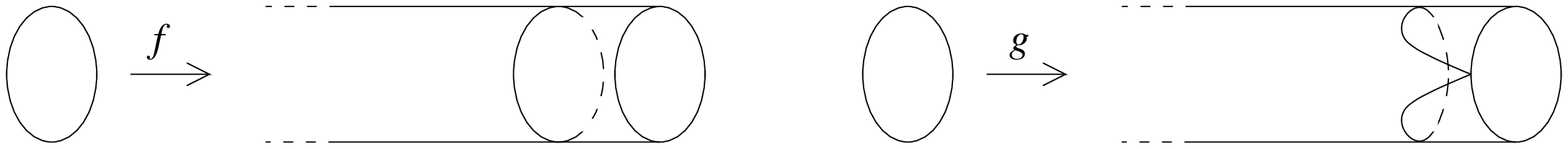,height=1cm}}
\vspace{.6cm}

They are both homotopic to $h(z)=(z,0)$. Now let $TM$ be
the tangent bundle of $M$ (it's a vector bundle in an orbifold sense
\footnote{We will not define orbifold vector bundles in this paper,
nor how they pull back. Therefore, we just hope that the facts
stated seem plausible.}).
The bundles $f^*(TM)$ and $g^*(TM)$ are respectively the orientable
and non-orientable 2-plane bundles
over $S^1$. By homotopy invariance 
of pullbacks, the bundle $h^*(TM)$ is isomorphic to both
$f^*(TM)$ and to $g^*(TM)$!

Of course, this just means that $h(z)=(z,0)$ is not a well defined
orbifold map.
\end{example}

Let us try carefully to understand the maps that are involved in this example.
The orbifold $M=S^1\times [\R/\Zt]$ is described by
$PM=S^1\times(\R\times E\Zt)/\Zt$, $QM=S^1\times\R/\Zt$ and
$p_M(z,[t,q])=(z,[t])$, where $[\:]$ means $\Zt$-orbit.

Choose a basepoint $\ell\in E\Zt$.
The map $f$ can be represented by
$$
\begin{CD}
S^1@>z\mapsto(z,[1,\ell])>Pf>PM&=S^1\times(\R\times E&\Z)/\Zt\\
@VVIdV@VVp_MV&\\
S^1@>z\mapsto(z,[1])>Qf>QM&=S^1\times\R/\Zt.&
\end{CD}
$$
Of all maps of example {\bf\ref{paradoxe}}, it's the only one that
was unambiguously defined. Indeed, all $p_M$-fibers over points of the form
$(z,[1])$ are contractible. Hence all orbifold maps $f':S^1\to M$
with $Qf'=Qf$ are homotopic to $f$ over $QM$ and are thus essentially the same
as $f$.

Choose a map $\gamma:S^1\setminus\{1\}\to E\Zt$ such that
$\lim_{\theta\to 0^+}\gamma(e^{i\theta})=\ell$ and
$\lim_{\theta\to 0^-}\gamma(e^{i\theta})=\overline{1}\cdot\ell$.
The orbifold map $g$ can be described by
$$
\begin{CD}
S^1@>\scriptstyle z\mapsto(z,[Re(\sqrt{z}),\gamma(z)])
\atop
\scriptstyle 1\mapsto(1,[1,\ell])\hspace{1.21cm}
>Pg>PM&=S^1\times(\R\times E&\Z)/\Zt\\
@VVIdV@VVp_MV&\\
S^1@>z\mapsto(z,[Re(\sqrt{z})])>Qg>QM&=S^1\times\R/\Zt,&
\end{CD}
$$
where we used the branch of $\sqrt{z}$ defined for $z\not\in\R_+$ and
satisfying $\sqrt{-1}=i$. One must
be careful not to confuse $g$ with the map $g'$ given by
$$
\begin{CD}
S^1@>z\mapsto(z,[|Re(\sqrt{z})|,\ell])
>Pg'>PM&=S^1\times(\R\times E&\Z)/\Zt\\
@VVIdV@VVp_MV&\\
S^1@>z\mapsto(z,[Re(\sqrt{z})])>Qg'>QM&=S^1\times\R/\Zt.&
\end{CD}
$$
The orbifold map $g$ is smooth (see definition {\bf\ref{orbifolddef}})
whereas $g'$ is not. In fact, every smooth orbifold map $g'':S^1\to M$
with $Qg''=Qg$ is homotopic to $Qg$ over $QM$.

Now, $f$ is homotopic to $h_0:S^1\to M$ given by
$$
\begin{CD}
S^1@>z\mapsto(z,[0,\ell])>Ph_0>PM&=S^1\times(\R\times E&\Z)/\Zt\\
@VVIdV@VVp_MV&\\
S^1@>z\mapsto(z,[0])>Qh_0>QM&=S^1\times\R/\Zt,&
\end{CD}
$$
whereas $g$ is homotopic to $h_1:S^1\to M$ given by
$$
\begin{CD}
S^1@>\scriptstyle z\mapsto(z,[0,\gamma(z)])\atop
\scriptstyle 1\mapsto(1,[0,\ell])\hspace{.4cm}
>Ph_1>PM&=S^1\times(\R\times E&\Z)/\Zt\\
@VVIdV@VVp_MV&\\
S^1@>z\mapsto(z,[0])>Qh_1>QM&=S^1\times\R/\Zt.&
\end{CD}
$$
The two orbifold maps $h_0$ and $h_1$ are not homotopy equivalent over
$QM$. To see that, consider the singular stratum of $M$. Call it $N$.
Its orbifold structure is described by $QN=S^1$, $PN=
S^1\times K(\mathbb Z_2,1)$, and $p_N$
the natural projection.
In this example we are interested in orbifold maps $h:S^1\to N$
that induce the identity on
the underlying spaces. These correspond to commutative diagrams
$$
\begin{CD}
S^1	@>Ph>>	S^1\times K(\mathbb Z_2,1)\\
@VVIdV	@VVpr_1V\\
S^1	@>Id=Qh>>	\phantom{,}S^1,
\end{CD}
$$
where $Ph$ is only defined up to homotopy over $S^1$.
Having such a diagram is the same as
having a section of $S^1\times K(\mathbb Z_2,1)\to S^1$ and these
are classified up to homotopy by
$$
[S^1,K(\mathbb Z_2,1)]=\mathbb Z_2.
$$
Therefore, we have essentially two orbifold maps
from $S^1$ to $N$ that induce the identity on the underlying spaces;
$h_0$ is one of them and
$h_1$ is the other.\medskip

Consider now the following example. It illustrates another subtlety
of the theory that escaped from the attention of many people.

\begin{example}
Let $S^1$ act on $S^3$ by $z\cdot(u,v)=(z^2u,z^2v)$ and $\mathbb Z_2$
act trivially on $S^2$. The two resulting orbifolds $[S^3/S^1]$ and
$[S^2/\mathbb Z_2]$ are both $S^2$'s with $\mathbb Z_2$ stabilizer
at every point. But they are not isomorphic.
\end{example}

Here is a way to see it. Every presentation of an orbifold $M$
as a global quotient $[X/G]$ gives a fibration $G\to X\times EG\to
(X\times EG)/G=PM$ and also, an associated long exact sequence of homotopy
groups. Defining $\pi_k(M)$ as $\pi_k(PM)$, this sequence then
looks like
$$
\ldots\pi_k(G)\to\pi_k(X)\to\pi_k(M)\to\pi_{k-1}(G)\ldots
$$
In our two examples the sequences we get are
$$
\ldots\underbrace{\pi_1(S^3)}_0\to\pi_1([S^3/S^1])\to\underbrace{\pi_0(S^1)}_0
\ldots
$$
and
$$
\ldots\underbrace{\pi_1(S^2)}_0\to\pi_1([S^2/\mathbb Z_2])\to
\underbrace{\pi_0(\mathbb Z_2)}_{\mathbb Z_2}\to
\underbrace{\pi_0(S^2)}_{0}\ldots
$$
This tells us that $\pi_1([S^3/S^1])=0$ and $\pi_1([S^2/\mathbb Z_2])
=\mathbb Z_2$ and in particular that these two orbifolds are not isomorphic.
With our definitions, the explanation is very easy, namely that the fibrations
$K(\mathbb Z_2,1)\to P[S^3/S^1]\to Q[S^3/S^1]$ and
$K(\mathbb Z_2,1)\to P[S^2/\mathbb Z_2]\to Q[S^2/\mathbb Z_2]$
are not isomorphic.

Note that it is very important to distinguish
them if one wants to have a good definition of homotopy groups of orbifolds.
It turns out that these two examples are
very good for testing one's definition of orbifolds.

\section {Definitions}

Let $p:P\to Q$ be a continuous map.

Given a point $x\in Q$, we say that a neighborhood
$\U\ni x$ is {\it good} if the inclusion $p^{-1}(x)\to p^{-1}(\U)$ is a
homotopy equivalence. An open set $\U$ is {\it good} if it is a good
neighborhood of some $x\in \U$.

We say that two maps $f,f':X\to P$ with $p\circ f=p\circ f'$ are
{\it homotopic over} $Q$ if there is a continuous family of maps
$f_t$, $t\in [0,1]$ with $f_0=f$, $f_1=f'$ and such that the composition
$p\circ f_t$ does not vary with $t$. We denote it by $f\simeq_Q f'$.

Given a map $\varphi:A\to B$, let us
denote by $\fp(A;B)$ the quotient $A/\!\!\sim$,
where $a\sim a'$ if $\varphi(a)=\varphi(a')$ and
if $a$ and $a'$ lie in the same
connected component of $\varphi^{-1}(\varphi(a))$.
The choice of $\varphi$ will always be cear by the context.

Let $\widetilde{X}$ stand for the universal cover of $X$.

\begin{df}\label{orbidef}
An \underline{orbispace} $M$ is a map $p_M:PM\to QM$ such that all $x\in QM$
have a basis of good neighborhoods $\U\ni x$ with the property that the maps of
universal covers
\begin{equation}\label{fibration}
\widetilde{p_M^{-1}(\U)}\quad
\longrightarrow\quad\fp\left(\widetilde{p_M^{-1}(\U)}\,;\,\U\right)
\end{equation}
are fibrations with contractible fibers
\footnote{The space
$\fp\big(\widetilde{p_M^{-1}(\U)};\U\big)$ should be thought as the local
branched cover over the open set $\U$.}.
The space $PM$ is called
the \underline{Borel space} and $QM$ the
\underline{underlying space} of $M$.

An \underline{orbispace map}
$f:M\to N$ is a map $Qf:QM\to QN$ and an equivalence class of maps
$Pf:PM\to PN$ that make the diagrams
$$
\begin{CD}
PM	@>Pf>>	PN\\
@VVp_MV	@VVp_NV\\
QM	@>Qf>>	QN
\end{CD}
$$
commute, under the equivalence relation of being homotopic over $QN$
\footnote{
The good way of formalizing the theory is to say that orbifolds form a
{\it 2-category}. This means that we don't consider two morphisms
$f=(Qf,Pf)$ and $g=(Qg,Pg)$ with $Qf=Qg$ and $Pf\simeq_{QN}Pg$ as {\it equal}
but we say instead that there's a {\it 2-morphism} $h:f\to g$. The 2-morphism
$h$ is the homotopy over $QN$.
}.


We say that two orbispaces $M$ and $N$
are \underline{isomorphic} if there exist orbispace maps $f:M\to N$
and $g:N\to M$ with $f\circ g = Id_N$
and $g\circ f = Id_M$. Note that $PM$ doesn't need to be homeomorphic to
$PN$ but only homotopy equivalent.

Given $x\in QM$ we call $\pi_1\big(p_M^{-1}(x)\big)$ the
\underline{stabilizer} of $x$ and denote it by $G_x$.
\end{df}

As it turns out, if $M$ is an orbispace and $\U$ {\it any} good open set,
then the map (\ref{fibration}) is a fibration with contractible fibers.

\begin{xx}[spaces]\rm
Any (locally contractible) space $X$, produces
an orbispace $Id_X:X\to X$.
Conversely, any orbispace $M$ with contractible fibers is isomorphic
to $Id_{QM}:QM\to QM$. We will say that such an orbispace ``is a space''. 
\end{xx}

\begin{prop}
Let $M$ be an orbispace. Then for all $x\in QM$, the fiber $p_M^{-1}(x)$ is 
a $K(G_x,1)$.
\end{prop}

\begin{proof} Let $\U\ni x$ be a good neighborhood of $x$. By definition
of good neighborhood, the homomorphism 
$G_x\equiv \pi_1\big(p_M^{-1}(x)\big)\to
\pi_1\big(p_M^{-1}(\U)\big)$ is an isomorphism.
Thus the fiber over $x$ of the map
$\widetilde{p_M^{-1}(\U)}\to\U$ is exactly $\widetilde{p_M^{-1}(x)}$
and in particular is connected. Let $\bar x$ be the corresponding point
in $\fp\Big(\widetilde{p_M^{-1}(\U)}\,;\,\U\Big)$. The space
$\widetilde{p_M^{-1}(x)}$ is also the fiber over $\bar x$ of the canonical map 
$$
\widetilde{p_M^{-1}(\U)}\quad
\longrightarrow\quad\fp\left(\widetilde{p_M^{-1}(\U)}\,;\,\U\right).
$$
It is contractible by definition of orbispaces hence
${p_M^{-1}(x)}$ is a $K(G_x,1)$.
\end{proof}

\begin{prop}
[global quotient]\label{global}
Let $G$ be a discrete group acting properly discontinuously on a
(locally equivariantly contractible\footnote{
By ``locally equivariantly contractible'', we mean that every point
$x$ in $X$ has a basis of $G_x$-invariant neighborhoods $\bar\U$.
Moreover the deformation retraction $\bar\U\searrow\{x\}$
should be realizable through $G_x$-equivariant maps. It's a
harmless assumption that is realized in all situations of interest:
simplicial actions on simplicial complexes, smooth actions
on manifolds...
}) space $X$. Then
$$
[X/G]\;=\;\bigg(\;p:(X\times EG)\big/G\to X/G\;\bigg)
$$
is an orbispace.
\end{prop}
\vspace{.6cm}

\centerline{\begin{minipage}{10cm}
\centerline{\epsfig{file=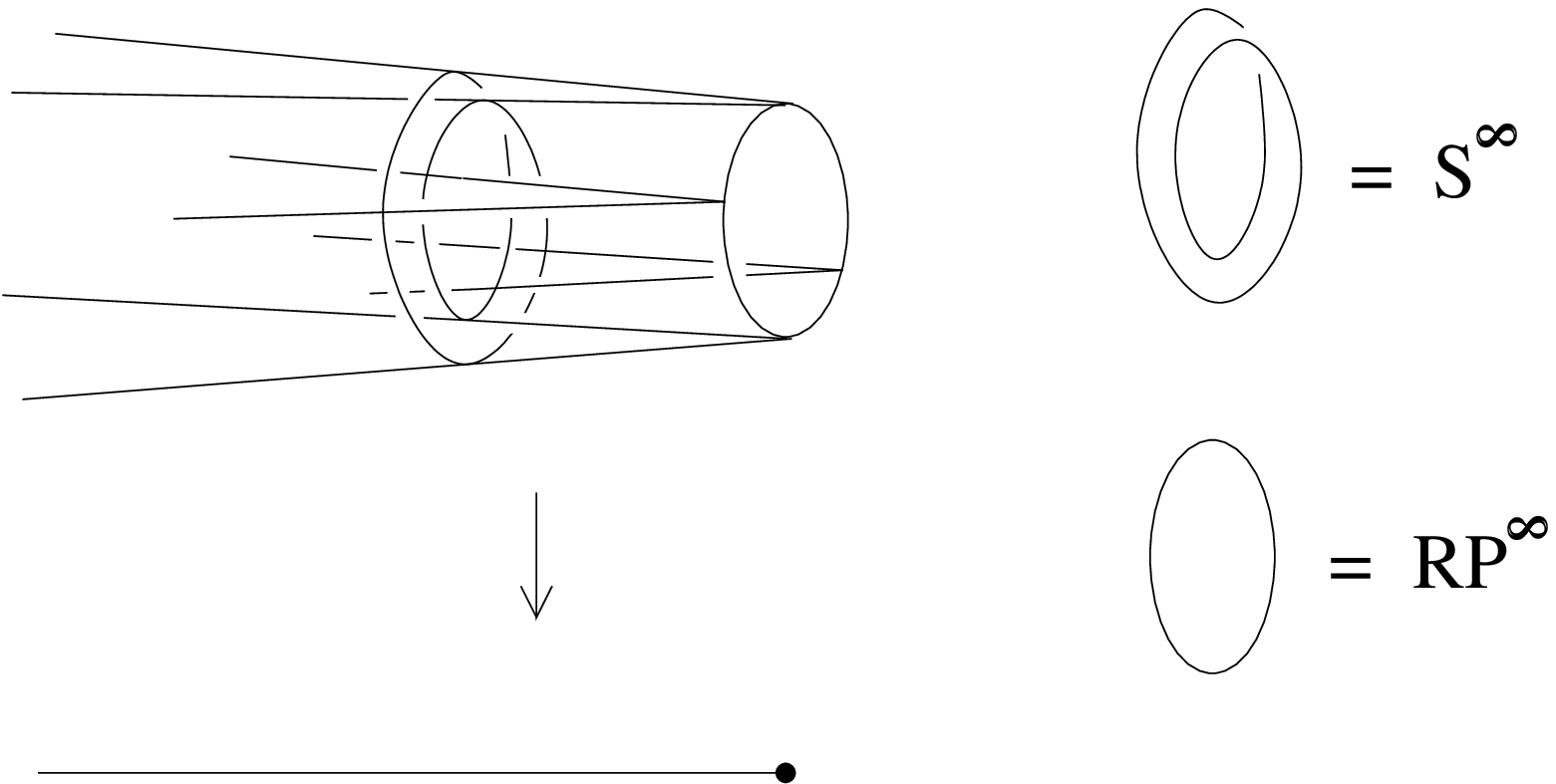,height=4cm}}
\vspace{.3cm}
\centerline{The orbispace $[\mathbb R/\mathbb Z_2]$ and its fibers.}
\vspace{.6cm}
\end{minipage}}

Note that different choices of $EG$ give rise to isomorphic orbispaces. 

\begin{proof} Let $\bar x\in X$ be a point, $x$ its image in $X/G$, and
$G_x$ its stabilizer under the action of G.
Note that $G_x=\pi_1\big(p^{-1}(x)\big)$ and 
that the two uses of the word {\it stabilizer} are consistent with
each other.

Let $\bar\U\ni \bar x$ be an equivariantly
contractible $G_x$-invariant neighborhood
that is small enough so that it doesn't meet it's
translates by $G$. Call $\U$ its image in
$X/G$. We claim that $\U$ is a neighborhood of $x$ with the
desired properties.
First we check that
\[
\begin{split}
p^{-1}(x)&
=(G\bar x\times EG)/G=(\{\bar x\}\times EG)/G_x\\
&\simeq(\bar\U\times EG)/G_x=(G\bar\U\times EG)/G=p^{-1}(\U).
\end{split}
\]
Hence $\U$ is a good neighborhood of $x$. Now
$$
\widetilde{p^{-1}(\U)}=\widetilde{(\bar\U\times EG)/G_x}=\bar\U\times EG
$$
and $\fp(\bar\U\times EG\,;\,\U)=\bar\U$, the map $\bar\U\times EG\to\bar\U$ is
indeed a fibration with contractible fibers.
We get a basis of good neighborhoods by repeating the
construction with smaller and smaller $\bar\U$.
\end{proof}\medskip

Any orbispace is locally isomorphic (in the sense of
definition {\bf \ref{orbidef}}) to a quotient constructed as in example
{\bf \ref{gl}}. To show this we need the following lemma.

\begin{lemma}
Let $G\acts B$ be any continuous action of a discrete group.
For $i\in\{1,2\}$, let
$G\acts E_i$ be free actions on $CW$-complexes, and let $p_i:E_i\to B$ be
equivariant fibrations with contractible fibers. Then these
fibrations are equivariantly homotopy equivalent over $B$.

More precisely, there are maps
$$
f:E_1\to E_2,\quad g:E_2\to E_1,\quad h_i:E_i\times[0,1]\to E_i
$$
such that all these maps are equivariant and that
$$
\begin{array}{ccc}
p_2\circ f=p_1,&\quad p_1\circ g=p_2,&
\quad p_i\circ h_i=p_i\circ{\rm pr}_1,\vspace{.2cm}\\
h_1\big|_{E_1\times\{0\}}=g\circ f,&\quad h_2\big|_{E_2\times\{0\}}=f\circ g,&
\quad h_i\big|_{E_i\times\{1\}}=Id_{E_i}.
\end{array}
$$\label{CW}
\end{lemma}

\begin{proof}[Idea of the proof]
Assume that $f$ is defined on the $n$-skeleton of $E_1$. Consider the
free $G$-action on the set of $(n+1)$-cells and choose representatives
of the orbits. One first extends the map $f$ to these cells and then
to all other cells by $G$-equivariance. This
process defines $f$, and the same for $g$.

Similarly one extends
$$
h_1\big|_{E_1\times\{0,1\}}:=g\circ f\sqcup Id_{E_1}:E_1\times\{0,1\}\to E_1
$$
to $h_1:E_1\times[0,1]\to E_1$,
and the same for $h_2$.
\end{proof}

To avoid technicalities, we will assume that all the spaces we deal with
are $CW$-complexes.
Thereafter we denote the local branched cover
$\fp\Big(\widetilde{p_M^{-1}(\U)}\,;\,\U\Big)$ by $\tU$.

\begin{prop}[local form for orbispaces]\label{loc1}
Let $M$ be an orbispace 
and $x\in QM$ a point.
Then there is a neighborhood
$\U\ni x$ such that the orbispace
$p_M^{-1}(\U)\to\U$ is isomorphic to
$\big[\hspace{.03cm}\tU/G_x\big]$.
\end{prop}

\begin{proof} Take $\U$ to be any good neighborhood
satisfying the requirements of definition {\bf\ref{orbidef}}.
Clearly the underlying space
$\tU/G_x$ is just $\U$.
The fibrations
$$
\tU\times EG_x\;\longrightarrow\;
\tU
$$
and
$$
\widetilde{p_M^{-1}(\U)}\;\longrightarrow\;
\tU
$$
satisfy the hypothesis of lemma {\bf\ref{CW}}.
Hence there are $G_x$-equivariant maps
$$
\tU\times EG_x
\hspace{.4cm}\overrightarrow{\longleftarrow}\hspace{.4cm}
\widetilde{p_M^{-1}(\U)}
$$
that are compatible with the projections on
$\tU$
and that are homotopy equivalences over $\U$.

Modding out by $G_x$, we get the desired orbispace maps
$$
\begin{array}{ccc}
\big(\tU\times EG_x\big)\big/G_x&
\overrightarrow{\longleftarrow}&p_M^{-1}(\U)\\
\downarrow&&\downarrow\vspace{.1cm}\\
\U&=\!\!=\!\!=&\U.
\end{array}
$$

\end{proof}

\begin{remark}
More generally, $[X/G]$ could have been defined to be
$X^+/G\to X/G$, where $X^+$ is a free $G-space$ and
$X^+\to X$ is an equivariant fibration with
contractible fibers. This would have made the statement of
proposition {\bf \ref{loc1}} almost tautological.
Note that lemma {\bf \ref{CW}} would still be needed in showing
that two different choices of $X^+$ give rise to
isomorphic orbispaces.
\end{remark}

Let us also classify orbispace maps locally. Given two actions
$G\acts X$, $H\acts Y$, a group homomorphism $\sigma:G\to H$ and a
$\sigma$-equivariant map $r:X\to Y$, call $[r/\sigma]:
[X/G]\to[Y/H]$ the orbispace map
$$
\begin{CD}
(X\times EG)/G@>>>(Y\times EH)/H\\
@VVV@VVV\\
X/G@>>>Y/H
\end{CD}
$$
induced by $r\times E\sigma:X\times EG\to Y\times EH$.

\begin{prop}[local form for maps]\label{loc2}
Let $f:M\to N$ be an orbispace map and $x\in QM$ a point. Call
$\sigma:G_x\to G_{Qf(x)}$ the group homomorphism induced by $f$.
Then there are two neighborhoods $\U\ni x$, $\V\ni Qf(x)$ and a
$\sigma$-equivariant map $\hat f:\tU\to\tV$ such that
$(Pf\big|_{p_M^{-1}(\U)},Qf\big|_\U)$ is isomorphic to $[\hat f/\sigma]$.
\end{prop}

\begin{proof} Choose good neighborhoods $\U$ and $\V$ such that $\U\subseteq
Qf^{-1}(\V)$. By proposition {\bf \ref{loc1}} we may assume
$M=[\tU/G]$ and $N=[\tV/H]$. The map $Pf:(\tU\times EG)/G\to
(\tV\times EH)/H$ induces $\sigma$ on the fundamental groups.
Now, taking universal covers gives us a $\sigma$-equivariant map
$\widetilde{Pf}:\tU\times EG\to\tV\times EH$ that descends to a map
$\hat f:\tU\to\tV$. It is clear by construction that
$f=[\hat f/\sigma]$.
\end{proof}\medskip

\begin{remark}\rm
Our definition of orbispaces is equivalent to the one Haefliger gave
in \cite{Hae90}, using topological groupoids (up to some issues
of local contractibility).
In his paper \cite{Hae84}, he shows how to construct a map $PM\to QM$
from a topological groupoid. On the other hand, given $PM\to QM$,
it is possible to recover the groupoid up to equivalence
\footnote{The topological groupoids considered here are such that the 
source and target maps from the space of morphisms
to the space of objects are local homeomorphisms.
Two topological groupoids $\mathcal G$ and $\mathcal G'$ are
(Morita-) equivalent if there is a third one $\mathcal G''$ and two
equivalences of categories $\mathcal G\leftarrow\mathcal G''\to\mathcal G'$
such that all maps involved are local homeomorphisms.}.

To do so, pick a cover of $QM$ by good open sets $\{\U_i\}$. For each $\U_i$
choose a section $s_i$ of (\ref{fibration}). Now let
$s:\bigsqcup_i\tU_i
\to PM$
be the composition of $\sqcup_i s_i$ with the natural projection on $PM$.
Given these choices, one defines a topological groupoid $\mathcal G$ by
\begin{alignat*}{1}
Ob(\mathcal G)\quad&:=\hspace{.3cm}\underset{i}{\bigsqcup}\;\;\tU_i
,\vspace{.3cm}\\
{\rm Hom}_{\mathcal G}(a,b)&:=\begin{cases}
\text{\begin{minipage}{4cm}\vspace{.11cm}homotopy classes of paths
$\gamma:[0,1]\to F_{s(a)}$ joining
$s(a)$ to $s(b)$\end{minipage}\hspace{.7cm}}&
\text{if $p_M(s(a))=p_M(s(b))$,}\vspace{.3cm}\\
\emptyset&\text{otherwise,}
\vspace{.1cm}\end{cases}
\end{alignat*}
where $F_z=p_M^{-1}(p_M(z))$ is the fiber of $z$. 
By working with refinements, it is not hard to show that different choices
of $\U_i$ and $s_i$ give equivalent groupoids.
\end{remark}

\begin{df}\label{orbifolddef}
An \underline{orbifold} $M$ is an orbispace that is locally isomorphic
to an orbispace of the form $[V/G]$ where $G\acts V$ is a linear
representation of a finite group.

More precisely for all $x\in QM$,
there is a good neighborhood $\U\ni x$, a representation $G\acts V$ and
an open subset $\V\subseteq V/G$, such that $p_M^{-1}(\U)\to\U$ is
isomorphic in the sense of definition {\bf \ref{orbidef}}
to $(\bar\V\times EG)/G\to \V$, where $\bar\V$ is the preimage of $\V$
in $V$.
\end{df}

Let $\{\U_i\}_{i\in I}$ be a cover of $QM$ by good open subsets and
$\varphi_i:\tU\to\Rm$ be homeomorphisms on their images. We say that
$\{\U_i,\varphi_i\}_{i\in I}$ form an {\it atlas} on $M$ if we have the
additional property that for all $i,j\in I$, $x\in \U_i\cap\U_j$,
one can find $k\in I$ such that $x\in\U_k\subseteq\U_i\cap\U_j$. 
An atlas is {\it smooth} if $\forall i,j$ with $\U_i\subseteq\U_j$, the
composition
$$
\Rn\supseteq\varphi_i(\tU_i)\overset{\varphi_i^{-1}}{\longrightarrow}
\tU_i\longrightarrow\tU_j\overset{\varphi_j}{\longrightarrow}\Rn
$$
is smooth, where $\tU_i\rightarrow\tU_j$ is any open inclusion
obtained by a lifting
of the natural map $\widetilde{p_M^{-1}(\U_i)}\to p_M^{-1}(\U_j)$
to a map $\widetilde{p_M^{-1}(\U_i)}\to \widetilde{p_M^{-1}(\U_j)}$.
Such a lifting always exists because $\widetilde{p_M^{-1}(\U_i)}$ is
contractible, but is in general not unique.
In the special case when $i=j$, this is equivalent to asking that
$$
\pi_1\big(p_M^{-1}(\U_i)\big)\:\acts\:
\tU_i\overset{\simeq}{\longrightarrow}\varphi_i(\tU_i)\subseteq\Rn
$$ 
is a smooth action.

\begin{df}
A \underline{smooth orbifold} is an orbifold with a chosen
maximal smooth atlas.

Let $(M,\{\U_i,\varphi_i\})$ and $(N,\{\V_j,\psi_j\})$ be
smooth orbifolds.
A \underline {smooth orbifold} \underline{map} $f:M\to N$
is an orbispace map  such
that for all $i,j$ with $Qf(\U_i)\subseteq\V_j$ the
composition of the maps
$$
\Rm\supseteq\varphi_i(\tU_i)\overset{\varphi_i^{-1}}{\longrightarrow}
\tU_i\longrightarrow\tV_j\overset{\psi_j}{\longrightarrow}\psi_j(\tV_j)
\subseteq\Rn
$$
is smooth, where $\tU_i\rightarrow\tV_j$ is any map obtained by a
lifting of
$\widetilde{p_M^{-1}(\U_i)}\to p_M^{-1}(\U_i)\overset{Pf}{\rightarrow}
p_N^{-1}(\V_j)$ to a map $\widetilde{p_M^{-1}(\U_i)}\to
\widetilde{p_N^{-1}(\V_j)}$.
\end{df}

As usual, a smooth atlas determines a unique maximal smooth atlas,
so it's enough to give {\it some} smooth atlas to define a smooth orbifold.

\begin{prop}[global quotient]
Let $G$ be a compact Lie group that acts quasi-freely on a smooth manifold $M$.
Then
$$
[X/G]\;=\;\bigg(\;p:(X\times EG)\big/G\to X/G\;\bigg)
$$
is a smooth orbifold.
\end{prop}

\begin{proof} Let $\bar x\in X$ be a point, $x$ its image in $X/G$, $\mathcal O_x$
its orbit and $G_x$ its stabilizer. Pick an invariant Riemannian metric
on $X$. Now let $\varepsilon>0$ be any number
smaller than the radius of injectivity of
the exponential map of the normal bundle $N\mathcal O_x\to X$.

Let $\bar\U$ be the image in $X$ of the
$\varepsilon$-ball of $N_x\mathcal O_x$,
and let $\U$ be its image in $X/G$.
By construction $g\bar\U=\bar\U\;\forall g\in G_x$ and
$g\bar\U\cap\bar\U=\emptyset\;\forall g\not\in G_x$. By
following the argument of example {\bf\ref{global}} we see that
$p^{-1}(\U)\to\U$ is equal to $(\bar\U\times EG)/G_x\to\bar\U/G_x$ and hence
$[X/G]$ is an orbifold (note that one can take $EG$ for $EG_x$).

The same argument also tells us that
\begin{alignat*}{1}
\tU\equiv\:&\fp\Big(\widetilde{p^{-1}(\U)}\,;\,\U\Big)=
\fp\Big(\widetilde{(G\bar\U\times EG)/G}\,;\,\U\Big)\\
&=\fp\Big(\widetilde{(\bar\U\times EG)/G_x}\,;\,\U\Big)=
\fp\Big(\bar\U\times EG\,;\,\U\Big)=\bar\U.
\end{alignat*}
The maps identifying the $\tU$'s with the $\bar\U$'s, composed with
any diffeomorphisms of the $\bar\U$'s to some open subset of $\Rn$
define a smooth atlas on $[X/G]$.
\end{proof}

\section{Conclusion}

The point of view on orbifolds we presented here will be further developed
in a second paper.
This first one was just meant as an invitation to translate
into this language all the notions that have been defined for
groupoids such as vector bundles, principal bundles...

The advantage of this approach is very clear for those interested in
homotopy properties of orbifolds. For example, the Van Kampen theorem
becomes completely trivial with our definition.

We also hope that this point of view will alow to formulate notions
that have not been introduced yet, for example that of group actions on
orbifolds.

\section{Acknowledgements}

I wish to thank my adviser Allen Knutson for his support and for the many
fruitful conversations that I had with him.

This paper was written while I was benefiting from a university grant offered
by B. Sturmfels and A. Knutson.

\end{document}